\documentclass[11pt, reqno]{amsart}

\usepackage[OT2, T1]{fontenc}

\usepackage[usenames,dvipsnames]{color}

\usepackage{url}
\usepackage{amsmath}
\usepackage{array}
\usepackage{graphicx}
\usepackage{amsfonts}
\usepackage{amssymb}
\usepackage{amsthm}
\definecolor{link}{RGB}{11,0,128}
\usepackage[colorlinks=true,citecolor=NavyBlue,linkcolor=Bittersweet,urlcolor=link]{hyperref}
\usepackage{paralist}
\usepackage{colonequals}	
\usepackage{color}
\usepackage{mathabx}		
\usepackage{amscd}
\usepackage[all,cmtip]{xy}	
\usepackage{sseq}			
\usepackage{verbatim}
\usepackage{parskip}
\usepackage{microtype}
\usepackage[in]{fullpage}
\usepackage{mathrsfs} 			
\usepackage{bm} 				
\usepackage{cleveref}
\usepackage{enumitem}
\usepackage{moreenum}			
\usepackage[alphabetic,lite]{amsrefs} 	
\usepackage{stmaryrd}	

\DeclareSymbolFont{cyrletters}{OT2}{wncyr}{m}{n}
\DeclareMathSymbol{\Sha}{\mathalpha}{cyrletters}{"58}

\newcommand{\gA}{\alpha}

\newcommand{\bF}{\mathbb{F}}
\newcommand{\bG}{\mathbb{G}}

\newcommand{\bQ}{\mathbb{Q}}

\newcommand{\bZ}{\mathbb{Z}}

\newcommand{\cA}{\mathcal{A}}
\newcommand{\cB}{\mathcal{B}}

\newcommand{\cG}{\mathcal{G}}
\newcommand{\cH}{\mathcal{H}}

\newcommand{\cN}{\mathcal{N}}
\newcommand{\cO}{\mathcal{O}}

\newcommand{\cQ}{\mathcal{Q}}

\newcommand{\fm}{\mathfrak{m}}
\newcommand{\fo}{\mathfrak{o}}

\newcommand{\sS}{\mathscr{S}}

\newcommand{\ra}{\rightarrow}

\newcommand{\xra}{\xrightarrow}

\newcommand{\hra}{\hookrightarrow}

\newcommand{\wt}{\widetilde}
\newcommand{\wh}{\widehat}

\newcommand{\pr}{^{\prime}}

\newcommand{\ce}{\colonequals}
\newcommand{\ov}{\overline}

\renewcommand{\b}{\textbf}
\newcommand{\surjects}{\twoheadrightarrow}

\newcommand{\tensor}{\otimes} 		

\renewcommand{\i}{^{-1}}
\renewcommand{\th}{^{\mathrm{th}}}

\providecommand{\p}[1]{\left(#1\right)}

\providecommand{\f}[2]{\frac{#1}{#2}}
\DeclareMathOperator{\Ker}{Ker}			
\DeclareMathOperator{\Coker}{Coker}		
\DeclareMathOperator{\im}{Im}			
\DeclareMathOperator{\Hom}{Hom}			
\DeclareMathOperator{\Char}{char}		
\DeclareMathOperator{\Ext}{Ext}			
\DeclareMathOperator{\loc}{loc}		%
\DeclareMathOperator{\Gal}{Gal}	
\DeclareMathOperator{\tr}{tr}		
\DeclareMathOperator{\inv}{inv}	
\DeclareMathOperator{\Res}{Res}		
\DeclareMathOperator{\Sel}{Sel}		
\DeclareMathOperator{\Pic}{Pic}		
\DeclareMathOperator{\Frob}{Frob}		
\DeclareMathOperator{\NS}{NS}		
\newcommand{\ba}{\begin{aligned}}
\newcommand{\ea}{\end{aligned}}
\newcommand{\be}{\begin{equation}}
\newcommand{\ee}{\end{equation}}
\newcommand{\pf}{\begin{proof}}
\newcommand{\bpf}{\begin{proof}}
\newcommand{\epf}{\end{proof}}
\newcommand{\bthm}{\begin{thm}}
\newcommand{\ethm}{\end{thm}}
\newcommand{\bprop}{\begin{prop}}
\newcommand{\eprop}{\end{prop}}
\newcommand{\bcor}{\begin{cor}}
\newcommand{\ecor}{\end{cor}}

\newcommand{\brem}{\begin{rem}}
\newcommand{\erem}{\end{rem}}
\newcommand{\brems}{\begin{rems} \hfill \begin{enumerate}[label=\b{\thesubsection.},ref=\thesubsection]}
\newcommand{\remi}{\addtocounter{subsection}{1} \item}
\newcommand{\erems}{\end{enumerate} \end{rems}}

\newcommand{\begs}{\begin{egs} \hfill \begin{enumerate}[label=\b{\thesubsection.},ref=\thesubsection]}
\newcommand{\egi}{\addtocounter{subsection}{1} \item}
\newcommand{\eegs}{\end{enumerate} \end{egs}}

\newcommand{\blem}{\begin{lemma}}
\newcommand{\elem}{\end{lemma}}
\newcommand{\bconj}{\begin{conj}}
\newcommand{\econj}{\end{conj}}
\newcommand{\bprob}{\begin{Problem}}
\newcommand{\eprob}{\end{Problem}}
\newcommand{\bq}{\begin{q}}
\newcommand{\eq}{\end{q}}
\newcommand{\benum}{\begin{enumerate}[label={(\alph*)}]}
\newcommand{\benuma}{\begin{enumerate}[label={(\arabic*)}]}
\newcommand{\benumr}{\begin{enumerate}[label={(\roman*)}]}
\newcommand{\eenum}{\end{enumerate}}
\newcommand{\bc}{}
\newcommand{\beg}{\begin{eg}}
\newcommand{\eeg}{\end{eg}}
\newcommand{\tst}{\textstyle}
\newcommand{\bcl}{\begin{claim}}
\newcommand{\ecl}{\end{claim}}
\newcommand{\lab}{\label}

\newcommand{\qq}{\quad\quad}
\newcommand{\qqq}{\quad\quad\quad}
\newcommand{\qqqq}{\quad\quad\quad\quad}

\theoremstyle{plain}
\newtheorem{thm}[subsection]{Theorem}
\Crefname{thm}{Theorem}{Theorems}

\Crefname{rethm}{Theorem}{Theorem}
\newtheorem{prop}[subsection]{Proposition}
\Crefname{prop}{Proposition}{Propositions}
\newtheorem{q}[subsection]{Question}
\Crefname{q}{Question}{Questions}
\Crefname{eg}{Example}{Examples}
\newtheorem{Problem}[subsection]{Problem}
\Crefname{Problem}{Problem}{Problems}
\newtheorem{conj}[subsection]{Conjecture}
\Crefname{conj}{Conjecture}{Conjectures}
\newtheorem{cor}[subsection]{Corollary}
\Crefname{cor}{Corollary}{Corollaries}

\newtheorem{lemma}[subsection]{Lemma}

\Crefname{subprop}{Proposition}{Propositions}
\Crefname{subcor}{Corollary}{Corollaries}
\Crefname{sublem}{Lemma}{Lemmas}

\theoremstyle{remark}
\newtheorem{claim}[equation]{Claim}
\Crefname{claim}{Claim}{Claims}

\Crefname{subrem}{Remark}{Remarks}

\theoremstyle{definition}
\newtheorem{eg}[subsection]{Example}
\newtheorem{rem}[subsection]{Remark}
\Crefname{rem}{Remark}{Remarks}
\newtheorem*{rems}{Remarks}
\newtheorem*{egs}{Examples}

\newtheoremstyle{subsection-tweak}
   {11pt}
   {3pt}%
   {}
   {}%
   {\bfseries}
   {}%
   {.5em}
   {\thmnumber{\@{#1}{}\@{#2}.}%
    \thmnote{~{\bfseries#3.}}}

\Crefname{innercustomconj}{Conjecture}{Conjecture}

\theoremstyle{subsection-tweak}
\newtheorem{pp}[subsection]{}
\newcommand{\bpp}{\begin{pp}}
\newcommand{\epp}{\end{pp}}

\numberwithin{equation}{subsection}

\usepackage[foot]{amsaddr}

\begin{document}
\author{K\k{e}stutis \v{C}esnavi\v{c}ius}
\title{$p$-Selmer growth in extensions of degree $p$}
\date{\today}
\address{Department of Mathematics, University of California, Berkeley, CA 94720-3840, USA}
\email{kestutis@berkeley.edu}

\begin{abstract} 
There is a known analogy between growth questions for class groups and for Selmer groups. If $p$ is a prime, then the $p$-torsion of the ideal class group grows unboundedly in $\bZ/p\bZ$-extensions of a fixed number field $K$, so one expects the same for the $p$-Selmer group of a nonzero abelian variety over $K$. This Selmer group analogue is known in special cases and we prove it in general, along with a version for arbitrary global fields.  
 \end{abstract}

\maketitle

\section{Introduction}

\bpp[Growth of class groups and of Selmer groups]
It is a classical theorem of Gauss that the $2$-torsion subgroup $\Pic(\cO_L)[2]$ of the ideal class group of a quadratic number field $L$ can be arbitrarily large. Although unboundedness of $\#\Pic(\cO_L)[p]$ for an odd prime $p$ is a seemingly inaccessible conjecture, \cite{BCHIS66}*{VII-12, Thm.~4} explains how to extend Gauss' methods to prove that $\#\Pic(\cO_L)[p]$ is unbounded if $L/\bQ$ ranges over the $\bZ/p\bZ$-extensions instead.

As explained in \cite{Ces15a}, growth questions for ideal class groups and for Selmer groups of abelian varieties are often analogous. It is therefore natural to hope that for a prime $p$ and a nonzero abelian variety $A$ over $\bQ$, the $p$-Selmer group $\Sel_p A_L$ can be arbitrarily large when $L/\bQ$ ranges over the $\bZ/p\bZ$-extensions. Our main result confirms this expectation.
\epp

\bthm[\Cref{p-Sel-unbdd}] \lab{main}
Let $p$ be a prime, $K$ a global field, and $A$ a nonzero abelian variety over $K$. If $A[p](\ov{K}) \neq 0$ (for instance, if $p \neq \Char K$) or if $A$ is supersingular, then
\[
\#\Sel_p A_{L}
\]
is unbounded when $L$ ranges over the $\bZ/p\bZ$-extensions of $K$.
\ethm

\brems 
\remi \lab{other-A}
In the excluded case when $A[p](\ov{K}) = 0$ and $A$ is not supersingular (when also $\Char K = p$ and $\dim A > 2$), there nevertheless is an $n \in \bZ_{> 0}$ that depends on $A$ such that
\[
\#\Sel_{p^n} A_L
\]
is unbounded when $L$ ranges over the $\bZ/p^n\bZ$-extensions of $K$, see \Cref{p-Sel-unbdd}.

\remi \lab{isog}
See \Cref{phi-unbdd} for a version of \Cref{main} for Selmer groups of arbitrary isogenies.

\remi \lab{class-grow}
The proof of \Cref{main} also reproves the unbounded growth of the $p$-torsion subgroup of the ideal class group in $\bZ/p\bZ$-extensions of $K$, see \Cref{CG-unbdd}.

\remi \lab{deeper}
If $K$ is a number field, then the analogue of \Cref{main} for $\bZ/n\bZ$-extensions with $1 < n < p$ seems to lie much deeper: by \cite{Ces15a}*{4.1 (b)}, this analogue would imply the conjectured unboundedness of $p$-torsion of ideal class groups of $\bZ/n\bZ$-extensions of a finite extension of $K$. If $\Char K > 0$, then \cite{Ces15a}*{5.5} proves an analogue of this type with $n = 2$ and $p \neq \Char K$.
\erems

\bpp[Previous investigations]
\Cref{main} was known in a number of special cases: 
\begin{itemize}
\item By Matsuno \cite{Mat09}*{4.5}, if $K = \bQ$ and $\dim A = 1$; 

\item
By Clark and Sharif \cite{CS10}*{Thm.~3}, if $p \neq \Char K$ and $\dim A = 1$; 

\item 
By Creutz \cite{Cre11}*{1.1} (which improves Clark \cite{Cla04}*{Thm.~7}), if $\Char K = 0$, Galois acts trivially on $\NS(A_{\ov{K}})$, and $A$ has a principal polarization coming from a $K$-rational divisor; 

\item 
By \cite{Ces15a}*{4.2}, if $p \neq \Char K$ and $A$ has $\bZ/p\bZ$ or $\mu_p$ as a $K$-subgroup, or if $A$ has everywhere semiabelian reduction and $\bZ/p\bZ$ as a $K$-subgroup. 
\end{itemize}

Unboundedness of $p$-Selmer sizes has also been observed in a number of other settings---typical variants include allowing arbitrary $L/K$ as long as $[L : K]$ is bounded and/or also varying $A$ as long as $\dim A$ is constant. For such results, see Cassels \cite{Cas64}, B\"{o}lling \cite{Bol75}, Kramer \cite{Kra83}, Fisher \cite{Fis01}*{Cor.~2}, Kloosterman and Schaefer \cite{KS03}*{Thm.~2}, Kloosterman \cite{Klo05}*{1.1}, Matsuno \cite{Mat07}*{5.1}, \cite{Mat09}*{Thm.~A}, and Bartel \cite{Bar10}*{1.1 and 4.4}. For an attempt to understand $p$-Selmer behavior in $\bZ/p\bZ$-extensions of number fields in the elliptic curve case, see Brau \cite{Bra14}. 
\epp

\bq
In \Cref{main}, is $\#\Sha(A_L)[p]$ also unbounded?
\eq

In their respective special cases of \Cref{main}, Clark--Sharif \cite{CS10} and Creutz \cite{Cre11} prove that the answer is `yes'.

\bpp[An overview of the proofs and of the paper] \lab{overview}
The proof of \Cref{main} is given in \S\ref{GSG} and is based on arithmetic duality, the key input being a general version of the Cassels--Poitou--Tate exact sequence. In \S\ref{ACPTS} we include a proof of this sequence that treats all global fields on an equal footing and circumvents well-known difficulties in positive characteristic by exploiting topologies carried by cohomology groups of local fields. The crucial topological input is closedness and discreteness of the image of a certain global-to-local pullback map. The analysis of this map in \S\ref{disc-image} rests in part on the results of \cite{Ces15d} and leads to several improvements to the literature on arithmetic duality in positive characteristic, notably to \cite{GA09}*{\S 4} and to \cite{Mil70a}. To be able to simultaneously prove the growth of class groups mentioned in Remark \ref{class-grow}, in \S\ref{FSG} we present a general framework for Selmer groups that extends the framework of Selmer structures of Mazur and Rubin to arbitrary global fields (in positive characteristic the unramified subgroups that play the decisive role in Selmer structures tend to be too small). For an overview of the auxiliary results isolated in \Cref{app-a,app-b}, see the introductions of the appendices.
\epp

\bpp[Notation] \lab{not}
The following notation is in place for the rest of the paper: 
\begin{itemize}
\item
$K$ is a global field; 

\item
If $\Char K = 0$, then $S$ is the spectrum of the ring of integers of $K$; 

\item
If $\Char K > 0$, then $S$ is the proper smooth curve over a finite field such that the function field of $S$ is $K$;

\item
A place of $K$ is denoted by $v$, the resulting completion by $K_v$, the ring of integers and the residue field of $K_v$ by $\cO_v$ and $\bF_v$, and the maximal ideal of $\cO_v$ by $\fm_v$;

\item
A place of a finite extension $K\pr$ of $K$ is denoted by $v\pr$.
\end{itemize}
As usual, $\mu_n$ denotes $\Ker(\bG_m \xra{n} \bG_m)$ and $\gA_p$ denotes the Frobenius kernel of the additive group $\bG_a$ over $\bF_p$.
For an $n \in \bZ_{\ge 0}$ and a $p$-divisible group (or a scheme) $X$ over a base $S$ of characteristic $p > 0$, we let $\Frob_{p^n,\, X/S}\colon X \ra X^{(p^n)}$ denote the $n$-fold relative Frobenius morphism of $X$ over $S$. Further notation is recorded in the beginning of each individual section.
\epp

As mentioned in \S\ref{overview}, topology carried by cohomology groups of local fields will play an important role in treating all $K$ at once. This topology is always taken to be the one defined in \cite{Ces15d}*{3.1--3.2}. To avoid cluttering the proofs with repetitive citations, in \S\ref{basic-topo} we gather the main topological properties that we will need. We will use these properties without explicit reference.

\bpp[$H^n(K_v, G)$ for a commutative finite $G$] \lab{basic-topo}
Fix a place $v$ of $K$, a commutative finite $K_v$-group scheme $G$, and an $n \in \bZ_{\ge 0}$. By \cite{Ces15d}*{3.5 (c), 3.6--3.8}, $H^n(K_v, G)$ is a locally compact Hausdorff abelian topological group that is discrete if $n \neq 1$. By \cite{Ces15d}*{3.5 (b)}, if $G$ is \'{e}tale (in particular, if $\Char K = 0$), then $H^1(K_v, G)$ is also discrete. By \cite{Ces15d}*{3.10}, if $v \nmid \infty$ and $\cG$ is a commutative finite flat $\cO_v$-model of $G$, then  the pullback map identifies $H^n(\cO_v, \cG)$ with a compact open subgroup of $H^n(K_v, G)$ (if $n \ge 2$, then $H^n(\cO_v, \cG) = 0$ by \cite{Toe11}*{3.4}). By \cite{Ces15d}*{4.2}, if $G$ fits into an exact sequence $0 \ra H \ra G \ra Q \ra 0$ of commutative finite $K_v$-group schemes, then the maps in the resulting cohomology sequence are continuous. 
\epp

\bpp[Conventions] \lab{conv}
We identify the nonarchimedean $v$ with the closed points of $S$. For a nonempty open $U \subset S$, writing $v\not \in U$ signifies that $v$ does not correspond to a closed point of $U$ (and hence could be archimedean). For a field $F$, a choice of its algebraic closure is denoted by $\ov{F}$. Fppf cohomology is denoted by $H^n$; identifications with \'{e}tale cohomology are implicit and use \cite{Gro68}*{11.7~1$^{\circ}$)}; further identifications with Galois cohomology are likewise implicit. If $v$ is archimedean, then we implicitly make the Tate modification: $H^n(K_v, -)$ is our shorthand for $\wh{H}^n(K_v, -)$ (this does not affect the $H^n$ with $n \ge 1$). Compactness of a topological space does not entail Hausdorffness, and `locally compact' means that every point has a compact neighborhood.
\epp

\subsection*{Acknowledgements}
I thank Julio Brau, Pete L. Clark, Brendan Creutz, Bjorn Poonen, and Yunqing Tang for helpful conversations or correspondence.


\section{Discreteness of the image of global cohomology} \lab{disc-image}


\bpp[The setup] \lab{setup}
Throughout \S\ref{disc-image}, we let $U$ be a nonempty open subscheme of $S$ and let $\cG$ be a commutative finite flat $U$-group scheme. The objects of study are the pullback map
\be \lab{loc-def}
\textstyle \loc^n(\cG)\colon H^n(U, \cG) \ra \bigoplus_{v\not\in U} H^n(K_v, \cG) \quad \text{for } n \in \bZ_{\ge 0}\quad \quad \text{ and its kernel}\quad \quad D^n(\cG).
\ee
We seek to show in \Cref{concl} that $\im(\loc^n(\cG))$ is closed and discrete and that $D^n(\cG)$ is finite.
\epp

\bprop \lab{Cart}
The following square is Cartesian (with injective maps as indicated):
\[
\xymatrix{
H^1(U, \cG) \ar@{^(->}[r] \ar[d] & H^1(K, \cG) \ar[d] \\
\prod_{v\in U} H^1(\cO_v, \cG) \ar@{^(->}[r] & \prod_{v\in U} H^1(K_v, \cG).
}
\]
\eprop

\bpf
This is a special case of \cite{Ces16a}*{4.3}.
\epf

\bprop \lab{im-closed}
The image $\im\p{\loc^n(\cG)}$ is closed for every $n \in \bZ_{ \ge 0}$.
\eprop

\bpf
Let $\cH$ be the Cartier dual of $\cG$. By \cite{Ces14d}*{5.3}, the images of $\loc^n(\cG)$ and $\loc^{2 - n}(\cH)$ are orthogonal complements under the sum of cup product pairings 
\[
H^n(K_v, \cG) \times H^{2 - n}(K_v, \cH) \ra H^2(K_v, \bG_m) \xra{\inv_v} \bQ/\bZ.
\]
By \Cref{cup-cont} and the discreteness of $H^2(K_v, \bG_m)$ supplied by \cite{Ces15d}*{3.5 (b)}, these pairings are continuous, so the claim follows.
\epf

\bprop\lab{neq-1}
The image $\im (\loc^n(\cG))$ is discrete for $n \neq 1$.
\eprop

\bpf
Even the target of $\loc^n(\cG)$ is discrete for $n \neq 1$.
\epf

\blem \lab{shrink}
Let $U\pr \subset U$ be a nonempty open.
\benum
\item \lab{shrink-a}
If $\im(\loc^1(\cG_{U\pr}))$ is discrete, then so is $\im(\loc^1(\cG))$.

\item \lab{shrink-b}
If $\im(\loc^1(\cG_{U\pr}))$ is discrete and $D^1(\cG_{U\pr})$ is finite, then $D^1(\cG)$ is finite.
\eenum
\elem

\bpf
Let $W$ be a compact neighborhood of $0$ in $\bigoplus_{v\not\in U} H^1(K_v, \cG)$, so 
\[
\tst W\pr \ce W \times \prod_{v\in U \setminus U\pr} H^1(\cO_v, \cG)
\]
is a compact neighborhood of $0$ in $\bigoplus_{v\not\in U\pr} H^1(K_v, \cG)$. By \Cref{im-closed}, the discrete $\im(\loc^1(\cG_{U\pr}))$ is closed in $\bigoplus_{v\not\in U\pr} H^1(K_v, \cG)$, so $\im(\loc^1(\cG_{U\pr})) \cap W\pr$ is finite. Thus, $\im(\loc^1(\cG)) \cap W$ is finite, too, $W \setminus \p{\im(\loc^1(\cG)) \cap \p{ W \setminus \{ 0\}}}$ exhibits $0$ as an isolated point of $\im(\loc^1(\cG))$, and \ref{shrink-a} follows.

For \ref{shrink-b}, \Cref{Cart} gives the inclusion $D^1(\cG_{U\pr}) \subset D^1(\cG)$ in $H^1(U\pr, \cG)$, whereas 
\[
[D^1(\cG) : D^1(\cG_{U\pr})] \le \#\p{\im(\loc^1(\cG_{U\pr})) \cap W\pr}. \qedhere
\]
\epf

\blem \lab{K-Kpr}
Let $K\pr/K$ be a finite separable extension and $U\pr$ the normalization of $U$ in $K\pr$. 
\benum
\item
If $\im(\loc^1(\cG_{U\pr}))$ is discrete, then so is $\im(\loc^1(\cG))$.

\item
If $D^1(\cG_{U\pr})$ is finite, then so is $D^1(\cG)$. 
\eenum
\elem

\bpf
Let $F\pr/F$ be either $K\pr/K$ or $K\pr_{v\pr}/K_v$ for places $v\pr\mid v$ of $K\pr$ and $K$. The kernel of the restriction 
\[
r \colon H^1(F, \cG) \ra H^1(F\pr, \cG)
\]
is finite, as one sees by using the separability to enlarge $F\pr/F$ to a finite Galois extension, and then identifying $\Ker r$ with $H^1(\Gal(F\pr/F), \cG(F\pr))$, which is finite by inspection.
\benum
\item
Let $W\pr$ be a compact neighborhood of $0$ in $\bigoplus_{v\pr\not\in U\pr} H^1(K_{v\pr}\pr, \cG)$. By \cite{Ces15d}*{2.7 (viii)}, the restriction $H^1(K_v, \cG) \ra H^1(K\pr_{v\pr}, \cG)$ is continuous for each $v\pr\mid v$, so there is a compact neighborhood $W$ of $0$ in $\bigoplus_{v \not\in U} H^1(K_v, \cG)$ lying in the preimage of $W\pr$. As in the proof of \Cref{shrink}, $\im(\loc^1(\cG_{U\pr})) \cap W\pr$ is finite and it suffices to prove the finiteness of $\im(\loc^1(\cG)) \cap W$, which follows by using in addition the finiteness of 
\[
\tst \qq \Ker(\bigoplus_{v \not\in U} H^1(K_v, \cG) \ra \bigoplus_{v\pr\not\in U\pr} H^1(K_{v\pr}\pr, \cG)).
\]

\item
The finite $\Ker(H^1(K, \cG) \xra{r} H^1(K\pr, \cG))$ contains $\Ker(D^1(\cG) \ra D^1(\cG_{U\pr}))$.
\qedhere
\eenum
\epf

\blem \lab{ext}
Suppose that $\cG$ fits into an exact sequence 
\[
0 \ra \cH \ra \cG \ra \cQ \ra 0
\]
of commutative finite flat $U$-group schemes.
\benum
\item
If $\im(\loc^1(\cH))$ and $\im(\loc^1(\cQ))$ are discrete and $D^1(\cQ)$ is finite, then $\im(\loc^1(\cG))$ is discrete.

\item
If $D^1(\cH)$ and $D^1(\cQ)$ are finite, then so is $D^1(\cG)$.
\eenum
\elem

\bpf
\hfill
\benum
\item
Let $W_\cQ$ be a compact neighborhood of $0$ in $\bigoplus_{v\not\in U} H^1(K_v, \cQ)$. By the discreteness of $\im \p{\loc^1(\cQ)}$ and \Cref{im-closed}, 
\[
\tst \qq \#\p{\im \p{\loc^1(\cQ)} \cap W_\cQ} < \infty.
\]
We combine this with the Hausdorffness of $\bigoplus_{v\not\in U} H^1(K_v, \cQ)$ to shrink $W_\cQ$ to ensure that 
\[
\qq \im \p{\loc^1(\cQ)} \cap W_\cQ = \{ 0 \}.
\]
We then let $W_\cG$ be a compact neighborhood of $0$ in the preimage of $W_{\cQ}$ in $\bigoplus_{v\not\in U} H^1(K_v, \cG)$.

\bcl \lab{WH-cpct}
The preimage $W_{\cH}$ of $W_{\cG}$ in $\bigoplus_{v\not\in U} H^1(K_v, \cH)$ is a compact neighborhood of $0$.
\ecl

\bpf
Due to continuity, $W_\cH$ is a closed neighborhood of $0$, so only its compactness requires proof. Each $x \in W_{\cH}$ has a compact neighborhood 
\[
\tst \qq W_x \subset  \bigoplus_{v\not\in U} H^1(K_v, \cH).
\]
By \cite{Ces15d}*{4.4~(c)~(3)}, the map
\[
\tst \qq \bigoplus_{v\not\in U} H^1(K_v, \cH) \ra \bigoplus_{v\not\in U} H^1(K_v, \cG)
\]
is closed, and hence open onto its image. The image of $W_\cH$ is a closed, and hence compact, subspace of $W_\cG$, so it is contained in the image of the union $Z$ of a finite number of the $W_x$. Since $Z$ is compact, so is the union of its $\p{\bigoplus_{v\not\in U} \cQ(K_v)}$-translates. This union contains the closed subset $W_\cH$, which is therefore~compact.
\epf

As in the proof of \Cref{shrink}, \Cref{im-closed} and \Cref{WH-cpct} ensure the finiteness of $\im\p{\loc^1(\cH)} \cap W_\cH$ and it suffices to prove the finiteness of $\im\p{\loc^1(\cG)} \cap W_\cG$. 

\bcl \lab{I-disc}
The image $I$ of $\im(\loc^1(\cH))$ in $\bigoplus_{v\not\in U} H^1(K_v, \cG)$ is closed and discrete.
\ecl

\bpf
For the closedness, it suffices to combine \Cref{im-closed} with loc.~cit. For the discreteness, it suffices to use the finiteness of $I \cap W_\cG$ inherited from $\im\p{\loc^1(\cH)} \cap W_\cH$.
\epf

\bcl \lab{LG-disc}
If $J \subset H^1(U, \cG)$ is the preimage of $D^1(\cQ)$, then $\loc^1(\cG)(J)$ is closed and~discrete.
\ecl

\bpf
Since 
\[
\qq [\loc^1(\cG)(J) : I] \le \#D^1(\cQ) < \infty,
\]
every subset of $\loc^1(\cG)(J)$ is a union of finitely many translates of subsets of $I$, and hence is closed in $\bigoplus_{v\not\in U} H^1(K_v, \cG)$ due to \Cref{I-disc}.
\epf

By construction, 
\[
\qq \im\p{\loc^1(\cG)} \cap W_\cG = \loc^1(\cG)(J) \cap W_\cG,
\]
so \Cref{LG-disc} gives the finiteness.

\item
Since 
\[
\qq [D^1(\cG) : D^1(\cG) \cap \im(H^1(U, \cH))] \le \#D^1(\cQ),
\]
it suffices to prove that the preimage $P$ of $D^1(\cG)$ in $H^1(U, \cH)$ is finite. For this, we use the inequality
\[
\tst [P : D^1(\cH)] \le \prod_{v\not\in U} \# \cQ(K_v). \qedhere
\]
\qedhere
\eenum
\epf

For use in the proof of \Cref{H1-disc}, we recall the following well-known lemma.

\blem \lab{ap-attack}
For a field $F$ and a commutative finite $F$-group scheme $G$, there is a finite separable extension $F\pr/F$ such that $G_{F\pr}$ is a successive extension of $F\pr$-group schemes that are isomorphic to $\bZ/m\bZ$ with $m \in \bZ_{> 0}$, or to $\mu_m$ with $m \in \bZ_{> 0}$, or to $\gA_p$ with $p = \Char F$ (where $\gA_0 \ce 0$).
\elem

\bpf
The claim is clear for \'{e}tale $G$, and hence, by passing to Cartier duals, also for $G$ of multiplicative type. Thus, the connected-\'{e}tale sequence allows us to assume that $G$ is connected and has a connected Cartier dual. By \cite{SGA3II}*{XVII, 4.2.1 ii) $\Leftrightarrow$ iv)}, such a $G$ is a successive extension of $\gA_p$'s.
\epf

\bthm \lab{H1-disc}
The image $\im(\loc^1(\cG))$ is discrete and the kernel $D^1(\cG)$ is finite.
\ethm

\bpf
\Cref{K-Kpr,ap-attack} reduce to the case when $\cG_K$ is a successive extension as in \Cref{ap-attack}. We spread out and use \Cref{shrink,ext} to reduce further to the cases of $\cG = \bZ/m\bZ$, of $\cG = \mu_m$, and of $\cG = \gA_p$, and we use these formulas to extend $\cG$ to a finite flat $S$-group scheme $\wt{\cG}$.

We set 
\[
\tst W \ce \bigoplus_{v\not\in U} H^1(\cO_v, \wt{\cG}) \qqq \text{ (with $H^1(\cO_v, \wt{\cG}) \ce H^1(K_v, \cG)$ for $v \mid \infty$),} 
\]
so $W$ is an open neighborhood of $0$ in $\bigoplus_{v\not\in U} H^1(K_v, \cG)$. By \Cref{Cart}, the preimage of $W$ in $H^1(U, \cG)$ is $H^1(S, \wt{\cG})$. Since $H^1(S, \wt{\cG})$ is finite, so are $D^1(\cG)$ and $\im(\loc^1(\cG)) \cap W$. 
\epf

\brems
\remi
The finiteness of $D^1(\cG)$ proved in \Cref{H1-disc} improves \cite{GA09}*{4.3}, which proved such finiteness after replacing $\cG$ by $\cG_{U\pr}$ for a sufficiently small nonempty open $U\pr \subset U$.

\remi
For the discreteness of $\im(\loc^1(\cG))$ to hold, not a single $v\not\in U$ can be omitted from the direct sum in the target of $\loc^1(\cG)$. For instance, for a prime $p$, the image of 
\[
\bF_p[t, t\i]/\bF_p[t, t\i]^p \cong H^1(\bF_p[t, t\i], \gA_p) \ra H^1(\bF_p((t)), \gA_p) \cong \bF_p((t))/\bF_p((t))^p
\]
 is dense rather than discrete.\footnote{The isomorphism $H^1(\bF_p((t)), \gA_p) \cong \bF_p((t))/\bF_p((t))^p$ is a homeomorphism by \cite{Ces15d}*{4.3 (b) and 4.5}.}
\erems

\begin{q}
Do the closedness and discreteness of the subset
\[
\tst \im(\loc^1(\cG)) \subset \bigoplus_{v\not\in U} H^1(K_v, \cG)
\]
continue to hold for a larger class of $U$-group schemes $\cG$ of finite type?
\end{q}

We turn to preliminaries needed for \Cref{Dn}.

\bpp[The cohomology with compact supports sequence] \lab{Hc}
We let $H^n_c$ denote the fppf cohomology with compact supports that takes into account the infinite primes, as defined in \cite{Mil06}*{III.0.6~(a)}. Loc.~cit.~provides the promised exact sequence, which reads
\be \lab{Hc-seq}
\textstyle \dotsb \ra H^n_c(U, \cG) \xra{x_{c}^n(\cG)} H^n(U, \cG) \xra{\loc^n(\cG)} \bigoplus_{v \not\in U} H^n(K_v, \cG) \xra{\delta_{c}^n(\cG)} H^{n + 1}_c(U, \cG) \ra \dotsb
\ee
and gives $D^n(\cG) = \im(x^n_c(\cG))$ and $\im(\loc^n(\cG)) = \Ker(\delta_c^n(\cG))$.
\epp

\bpp[Global duality] \lab{GDP}
We let $\cH$ be the Cartier dual of $\cG$, so that \cite{Mil06}*{III.3.2 and III.8.2} gives a bilinear pairing 
\be \lab{GDP-def}
H^n(U, \cG) \times H^{3 - n}_c(U, \cH) \ra H^3_c(U, \bG_m) \xra{\tr} \bQ/\bZ
\ee
that identifies $H^{3 - n}_c(U, \cH)$ with the Pontryagin dual of the discrete $H^n(U, \cG)$.
\epp

\blem \lab{map-id}
In the setup of \S\S\ref{Hc}--\ref{GDP}, the homomorphism dual to $\loc^n(\cG)$ identifies with $\delta_c^{2 - n}(\cH)$, i.e., the following diagram commutes
\[
\xymatrix{
H^n(U, \cG) \ar[d]_-{\loc^n(\cG)} \ar@{}[r]|-{\bigtimes} & H^{3- n}_c(U, \cH) \ar[rr]^{\eqref{GDP-def}} && H^3_c(U, \bG_m) \ar[rr]^-{\tr} &&\bQ/\bZ \ar@{=}[d] \\
\bigoplus_{v\not\in U} H^n(K_v, \cG) \ar@{}[r]|-{\bigtimes} & \bigoplus_{v\not\in U} H^{2 - n}(K_v, \cH) \ar[u]^-{\delta^{2 - n}_c(\cH)} \ar[rr]^-{\sum_v - \cup -} && \bigoplus_{v\not\in U} H^2(K_v, \bG_m) \ar[u]^-{\delta^2_c(\bG_m)} \ar[rr]^-{\sum_v \inv_v} && \bQ/\bZ,
}
\]
where the bottom row is the sum of Tate--Shatz local duality pairings.
\elem

\bpf
We have proved this in the course of the proof of \cite{Ces14d}*{5.3}.
\epf

\bthm \lab{Dn}
The pairing \eqref{GDP-def} induces a perfect pairing of finite abelian groups
\be \lab{Dn-pair}
D^n(\cG) \times D^{3 - n}(\cH) \ra \bQ/\bZ.
\ee
\ethm

\bpf
\Cref{map-id} proves that 
\[
\Ker(\loc^n(\cG)) \qq \text{and} \qq \Ker(x^{3 - n}_c(\cH))
\]
are orthogonal under \eqref{GDP-def}, so \eqref{Dn-pair} exists and is nondegenerate on the left. Consequently, 
\[
\#D^n(\cG) \le \#D^{3 - n}(\cH),
\]
and, since $D^0(\cH)$ inherits finiteness from $H^0(U, \cH)$ and $D^1(\cH)$ is finite by \Cref{H1-disc}, $D^2(\cG)$ and $D^3(\cG)$ are finite, too. Swapping the roles of $\cG$ and $\cH$, we learn that equalities must hold in the inequalities above, so \eqref{Dn-pair} is also nondegenerate on the right.
\epf

\brem
The deduction of \Cref{Dn} from \Cref{H1-disc} is the same as that of \cite{GA09}*{4.7} from \cite{GA09}*{4.3}.
\erem

For ease of reference, we combine some of the results of \S\ref{disc-image} into the following theorem.

\bthm \lab{concl}
For $U$ and $\cG$ as in \S\ref{setup} and every $n\in \bZ_{\ge 0}$, the image $\im(\loc^n(\cG))$ is closed and discrete in $\bigoplus_{v\not\in U} H^n(K_v, \cG)$ and the kernel $D^n(\cG)$ is finite.
\ethm

\bpf
The image claim is proved in \Cref{im-closed,neq-1} and \Cref{H1-disc}. The kernel claim is proved in \Cref{H1-disc,Dn}.
\epf


\section{Finiteness of Selmer groups} \lab{FSG}

The results of \S\ref{disc-image} allow us to prove finiteness of Selmer groups without distinguishing between the number field and the function field cases (see \Cref{Sel-fin}). The key finiteness inputs to \Cref{Sel-fin} are the finiteness of class groups and the Dirichlet Unit Theorem, both through the proof of \Cref{H1-disc}. In the case of the $\phi$-Selmer group for an isogeny $\phi$ between abelian varieties over $K$, \Cref{Sel-fin} seems to improve the literature by treating all $K$ and $\phi$ simultaneously, instead of resorting to \cite{Mil70a} that was tailored specifically to the $\Char K \mid \deg \phi$ case. The approach of loc.~cit.~is close to ours: the key lemma of \cite{Mil70a} is a variant of \Cref{H1-disc} for $\bZ/p\bZ$, $\mu_p$, and~$\gA_p$. 

Throughout \S\ref{FSG}, we fix a commutative finite $K$-group scheme $G$.

\bpp[Selmer groups] \lab{Sel-gp}
\emph{Selmer conditions} for $G$ are compact subgroups
\be \lab{sel-cond}
\Sel(G_{K_v}) \subset H^1(K_v, G), \quad\quad \text{ one for each place $v$ of $K$,}
\ee
such that there is a nonempty open $U \subset S$ and a commutative finite flat $U$-model $\cG$ of $G$ for which 
\[
\Sel(G_{K_v}) \subset H^1(\cO_v, \cG) \quad \text{inside} \quad H^1(K_v, G) \qqq \text{for every} \quad v\in U
\]
(the choice of $\cG$ plays no role: two $\cG$'s identify over a smaller $U$). The resulting \emph{Selmer group}, $\Sel(G)$, is the fiber product
\[
\xymatrix{
\Sel(G)  \ar[d] \ar@{}[r]|-{\text{\scalebox{1.7}{$\subset$}}} & H^1(K, G) \ar[d] \\
\prod_{v} \Sel(G_{K_v}) \ar@{}[r]|-{\text{\scalebox{1.7}{$\subset$}}} & \prod_{v} H^{1}(K_v, G).
}
\]
(Implicitly, $\Sel(G)$ depends on the chosen Selmer conditions \eqref{sel-cond}.)
\epp

\bthm \lab{Sel-fin}
For every choice of Selmer conditions, $\Sel(G)$ is finite.
\ethm

\bpf
Let $U$ and $\cG$ be as in \S\ref{Sel-gp}. By \Cref{Cart}, imposing Selmer conditions at all $v\in U$ leaves us with a subgroup of $H^1(U, \cG)$. By \Cref{concl} and the compactness of $\bigoplus_{v\not\in U} \Sel(G_{K_v})$, imposing the further conditions at all $v\not\in U$ leaves us with a finite group.
\epf

\begs
\egi
If $\Sel(G_{K_v}) = 0$ for all $v$, then \Cref{Sel-fin} recovers the finiteness of 
\[
\textstyle\Sha^1(G) \ce \Ker(H^1(K, G) \ra \prod_v H^1(K_v, G)),
\]
proved in \cite{Mil06}*{I.4.9} in the number field case and in \cite{GA09}*{4.6} in the function field case.

\egi \lab{eg-phi}
If $\phi\colon A \ra B$ is an isogeny of abelian varieties over $K$, then the subgroups 
\[
\qq B(K_v)/\phi A(K_v) \subset H^1(K_v, A[\phi])
\]
 are compact due to the compactness of $B(K_v)$ and the continuity of the connecting map (supplied by \cite{Ces15d}*{4.2}). These subgroups are Selmer conditions---the $U$-model requirement  is met due to \cite{Ces16a}*{2.5~(d)}. The resulting Selmer group is the \emph{$\phi$-Selmer group} $\Sel_\phi A$. 

\egi \lab{eg-ff}
If $\cG$ is a commutative finite flat $S$-group scheme, then \Cref{Cart} ensures that $H^1(S, \cG)$ is the Selmer group that results from the Selmer conditions 
\[
\ba
\qqq &H^1(\cO_v, \cG) \subset H^1(K_v, \cG) \qq \text{for $v\nmid \infty$} \quad \text{and}\\
&H^1(K_v, \cG) \subset H^1(K_v, \cG) \qq \text{for $v\mid \infty$.}
\ea
\]

\egi \lab{eg-cl}
If $\Char K = 0$ and in Example \ref{eg-ff} one chooses $\cG = \bZ/m\bZ$ for $m\in \bZ_{\ge 0}$ but alters the Selmer conditions to be $0 \subset H^1(K_v, \cG)$ for $v\mid \infty$, then, by the theory of the Hilbert class field, the resulting Selmer group is the Pontryagin dual of the $m$-torsion of the ideal class group of $K$.
\eegs

\brem
The notion of Selmer conditions extends the notion of a Selmer structure defined in \cite{MR07}*{1.2} in a number field setting. The role of the $U$-model $\cG$  is analogous to the role of the unramified subgroups in loc.~cit.~(the unramified subgroups are too small when $G$ is not \'{e}tale), with the caveat that for added flexibility we do not insist that almost all of the inclusions $\Sel(G_{K_v}) \subset H^1(\cO_v, \cG)$ be equalities.
\erem


\section{Cassels--Poitou--Tate} \lab{ACPTS}

In \S\ref{GSG}, our proof of unbounded Selmer growth is based on manipulating a generalization of the Cassels--Poitou--Tate sequence. This generalization is presented in \Cref{PT-gen}, which extends \cite{CS00}*{1.5} to finite group schemes over global fields (loc.~cit.~focused on the case of finite group schemes of odd order over number fields). In \eqref{Coates} we write out the sequence of \Cref{PT-gen} in the special case of Selmer groups of dual isogenies between abelian varieties over a global field.

\bpp[Selmer conditions that are orthogonal complements] \lab{OC}
Let $U\subset S$ be a nonempty open, $\cG$ a commutative finite flat $U$-group scheme, and $\cH$ its Cartier dual. For each $v\not \in U$, let
\be \lab{OC-cond}
\Sel(\cG_{K_v}) \subset H^1(K_v, \cG) \quad\quad\text{and}\quad\quad \Sel(\cH_{K_v}) \subset H^1(K_v, \cH) \quad \quad 
\ee
be compact subgroups that are orthogonal complements under the Tate--Shatz local duality pairing
\be\lab{Shatz}
H^1(K_v, \cG) \times H^{1}(K_v, \cH) \ra H^2(K_v, \bG_m) \hra \bQ/\bZ,
\ee
which is perfect by \cite{Sha64}*{Duality theorem on p.~411} (alternatively, by \cite{Mil06}*{I.2.3, I.2.13 (a), III.6.10}). We complete \eqref{OC-cond} to Selmer conditions by using the compact subgroups
\be \lab{comp-cond}
\quad \quad \quad H^1(\cO_v, \cG) \subset H^1(K_v, \cG) \quad\quad\text{and}\quad\quad H^1(\cO_v, \cH) \subset H^1(K_v, \cH) \quad \quad \quad \text{for $v\in U$}.
\ee
By \cite{Mil06}*{III.1.4 and III.7.2}, \eqref{comp-cond} also concerns orthogonal complements, so shrinking $U$ does not affect the setup. By \Cref{Cart}, the resulting Selmer groups $\Sel(\cG)$ and $\Sel(\cH)$ fit into inclusions
\be \lab{fit-into}
\Sel(\cG) \subset H^1(U, \cG) \subset H^1(K, \cG) \quad \quad \text{and} \quad \quad \Sel(\cH) \subset H^1(U, \cH) \subset H^1(K, \cH); \quad  \quad \quad \ \  
\ee
by \Cref{Sel-fin}, they are finite. As in \S\ref{setup}, we let 
\[
\tst \loc^n(\cG) \colon H^n(U, \cG) \ra \bigoplus_{v\not\in U} H^n(K_v, \cG)
\]
be the pullback map.
\epp

\bthm \lab{PT-gen} 
With the setup of \S\ref{OC} there is an exact sequence with continuous maps
\[
0 \ra \Sel(\cG) \ra H^1(U, \cG) \ra \bigoplus_{v\not\in U} \f{H^1(K_v, \cG)}{\Sel(\cG_{K_v})} \xra{y(\cG)} \Sel(\cH)^* \xra{x(\cG)} H^2(U, \cG) \xra{\loc^2(\cG)} \bigoplus_{v\not\in U} H^2(K_v, \cG),
\]
where $\Sel(\cG)$, $\Sel(\cH)$, $H^1(U, \cG)$, and $H^2(U, \cG)$ are discrete and $(-)^*$ denotes the Pontryagin dual.
\ethm

\bpf
Exactness up to $H^1(U, \cG)$ amounts to \eqref{fit-into} and the definition of $\Sel(\cG)$. By \cite{Ces14d}*{5.3}, 
\[
\textstyle \im(\loc^1(\cG)) \subset \bigoplus_{v\not\in U} H^1(K_v, \cG) \quad \text{ and } \quad \im(\loc^{1}(\cH)) \subset \bigoplus_{v\not\in U} H^1(K_v, \cH)
\]
 are (closed) orthogonal complements under the sum of the pairings \eqref{Shatz}, and hence, by \cite{BouTG}*{III.28, Cor.~1} and \cite{HR79}*{24.10}, so are
\[
\im(\loc^1(\cG))\,\textstyle{+} \bigoplus_{v\not\in U} \Sel(\cG_{K_v})  \subset \bigoplus_{v\not\in U} H^1(K_v, \cG)  \quad\text{and}\quad \im(\loc^1(\cH)|_{\Sel(\cH)}) \subset \bigoplus_{v\not\in U} H^1(K_v, \cH).
\]  
We therefore arrive at further orthogonal complements
\be\lab{so-are}\tag{$\dagger$}
\textstyle{}\im\p{H^1(U, \cG) \ra  \bigoplus_{v\not\in U} \f{H^1(K_v, \cG)}{\Sel(\cG_{K_v})} }  \quad\text{and}\quad \im(\loc^1(\cH)|_{\Sel(\cH)}) \subset \bigoplus_{v\not\in U} \Sel(\cH_{K_v}),
\ee
which are closed because \cite{HR79}*{24.11} ensures the continuity of the pairings between $\f{H^1(K_v, \cG)}{\Sel(\cG_{K_v})}$ and $\Sel(\cH_{K_v})$. Loc.~cit.~then allows us to define $y(\cG)$ to be the continuous map that factors through 
\[
\im(\loc^1(\cH)|_{\Sel(\cH)})^* \ra \Sel(\cH)^*.
\]
This map is injective due to \Cref{concl} and \cite{HR79}*{24.8--24.11}, so exactness at $\bigoplus_{v\not\in U} \f{H^1(K_v, \cG)}{\Sel(\cG_{K_v})}$ follows.

Loc.~cit.~also ensures the exactness of the sequence
\[
\tst 0 \ra \p{\f{H^1(U, \cH)}{\Sel(\cH)}}^* \ra H^1(U, \cH)^* \ra \Sel(\cH)^* \ra 0.
\]
By \Cref{concl} and \cite{BouTG}*{III.28, Cor.~3}, the image of 
\[
\tst H^1(U, \cH) \ra\bigoplus_{v\not\in U} \f{H^1(K_v, \cH)}{\Sel(\cH_{K_v})}
\]
is discrete. By the analogue of \eqref{so-are} for $\cH$, this image is also closed. Thus, by \cite{HR79}*{23.18, 24.8, and 24.11}, the middle row of the diagram
\[
\xymatrix{
&& \bigoplus_{v\not\in U} H^1(K_v, \cG) \ar@{->>}[r] \ar[d]^-{\loc^1(\cH)^*} & \bigoplus_{v\not\in U} \f{H^1(K_v, \cG)}{\Sel(\cG_{K_v})} \ar[d]^-{y(\cG)} \\
\textstyle \bigoplus_{v\not\in U} \Sel(\cG_{K_v}) \ar@{^(->}[rru] \ar@{}[r]|-{\text{\scalebox{1.4}{$\cong$}}}  & \bigoplus_{v\not\in U} \p{\f{H^1(K_v, \cH)}{\Sel(\cH_{K_v})}}^* \ar[r] & H^1(U, \cH)^* \ar[r] \ar[d]^-{x} & \Sel(\cH)^* \ar[r] \ar@{-->}[ld]^(0.45){x(\cG)} & 0 \\
&& H^2(U, \cG)
}
\]
is exact. Since $\loc^1(\cH)|_{\Sel(\cH)}$ factors through $\bigoplus_{v\not\in U} \Sel(\cH_{K_v})$, the top part of the diagram commutes. The map $x$ is obtained from the map $x^2_c(\cH)$ of \eqref{Hc-seq} by using \eqref{GDP-def}, so the middle column is exact by \Cref{map-id}. In conclusion, $x$ factors through a unique $x(\cG)$ as indicated, 
\[
\Ker(x(\cG)) = \im(y(\cG)),\qq  \text{and} \qq  \im(x(\cG)) = \im x = \Ker(\loc^2(\cG)).
\]
Finiteness of $\Sel(\cH)^*$ ensures the continuity of $x(\cG)$.
\epf

\brems
\remi \lab{right}
To extend the sequence of \Cref{PT-gen} to the right, combine \eqref{Hc-seq}, \S\ref{GDP}, and~\Cref{map-id}.

\remi \lab{over-S}
If $\cG$ and $\cH$ extend to Cartier dual finite flat $S$-group schemes $\wt{\cG}$ and $\wt{\cH}$, then one may take 
\[
\quad \quad \quad \Sel(\cG_{K_v}) = \begin{cases}H^1(\cO_v, \wt{\cG}), \quad  \text{for $v\in S \setminus U$,}  \\   H^1(K_v, \cG), \quad \text{for $v\mid \infty$,}\end{cases}  \text{and} \quad  \Sel(\cH_{K_v}) = \begin{cases}H^1(\cO_v, \wt{\cH}), \quad \text{for $v\in S \setminus U$,}  \\  0, \quad\quad\quad\quad\, \quad \text{for $v\mid \infty$.}\end{cases}
\]
With these choices, the sequence of \Cref{PT-gen} compares $H^1(S, \wt{\cG})$ and $H^1(U, \cG)$.
\erems

\beg \lab{CPT-phi}
Let $\phi\colon A \ra B$ and $\phi^\vee\colon B^\vee \ra A^\vee$ be dual isogenies of abelian varieties over $K$, and let $\phi\colon \cA \ra \cB$ and $\phi^\vee \colon \cB^\vee \ra \cA^\vee$ be their extensions to dual isogenies of abelian schemes over a nonempty open $U \subset S$ of good reduction. As in Example \ref{eg-phi}, 
\[
\quad \quad \quad B(K_v)/\phi A(K_v) \subset H^1(K_v, A[\phi]) \quad \quad \text{and} \quad \quad A^\vee(K_v)/\phi^\vee B^\vee(K_v) \subset H^1(K_v, B^\vee[\phi^\vee])\,
\]
constitute Selmer conditions for $A[\phi]$ and $B^\vee[\phi^\vee]$ with Selmer groups $\Sel_\phi A$ and $\Sel_{\phi^\vee} B^\vee$. By \Cref{Cart} and \cite{Ces16a}*{2.5~(d)}, these conditions for $v\in U$ cut out 
\[
\quad \quad \  H^1(U, \cA[\phi]) \subset H^1(K, A[\phi]) \quad \quad \text{and} \quad \quad H^1(U, \cB^\vee[\phi^\vee]) \subset H^1(K, B^\vee[\phi^\vee]).
\]
By \Cref{boring-check}, the remaining conditions at $v\not\in U$ put us in the framework \Cref{PT-gen}. Taking into account Remark \ref{right}, the resulting exact sequence reads
\be\ba \lab{Coates}
\xymatrix@C=12pt @R=8pt{
 0 \ar[r] & \Sel_{\phi} A \ar[r]  &  H^1(U, \cA[\phi]) \ar[r] &  \bigoplus_{v\not\in U} H^1(K_v, A)[\phi]  \ar[r] & (\Sel_{\phi^\vee} B^\vee)^* 
                \ar@{->} `r/7pt[d] `/8pt[l] `^dl[lll] `^r/6pt[dll] [dll] \\
             && H^2(U, \cA[\phi]) \ar[r] & \bigoplus_{v\not\in U} H^2(K_v, A[\phi]) \ar[r] & (B^\vee[\phi^\vee](K))^*. &
}
\ea\ee
The last map is surjective if $H^3(U, \cA[\phi]) = 0$. By \S\ref{GDP}, this is so if $H^0_c(U, \cB^\vee[\phi^\vee]) = 0$, in particular, if $U \neq S$ and either $\Char K = 0$ and $2 \nmid \deg \phi$, or $\Char K > 0$.
\eeg


\section{Growth of Selmer groups} \lab{GSG}

\Cref{Sel-gr,p-Sel-unbdd} along with \Cref{phi-unbdd} are the sought unbounded Selmer growth results. In \S\ref{GSG}, for a finite extension $K\pr/K$ and an open $U \subset S$, we denote the normalization of $U$ in $K\pr$ by $U\pr$.

\bpp[Selmer conditions over varying base fields] \lab{Sel-BC}
To fix the general setup, suppose that $U \subset S$ is a nonempty open, $\cG$ and $\cH$ are Cartier dual commutative finite flat $U$-group schemes, and $\sS$ is a set of finite extensions of $K$ such that $K \in \sS$. Suppose also that for each $K\pr \in \sS$ one has compact subgroups
\be \lab{Sel-Sel}
 \quad\quad \quad\quad  \ \,\Sel(\cG_{K\pr_{v\pr}}) \subset H^1({K\pr_{v\pr}}, \cG) \quad\quad\text{and}\quad\quad \Sel(\cH_{K\pr_{v\pr}}) \subset H^1({K\pr_{v\pr}}, \cH) \quad \quad \text{for $v\pr \not\in U\pr$}
\ee
that are orthogonal complements as in \eqref{OC-cond} and such that the restriction maps 
\[
H^1(K_v, \cG) \ra H^1({K\pr_{v\pr}}, \cG) \quad \quad \text{and} \quad\quad H^1(K_v, \cH) \ra H^1({K\pr_{v\pr}}, \cH) \quad \quad \quad \quad \ \ \,
\]
induce the maps
\be \lab{Sel-comp}
\Sel(\cG_{K_{v}}) \ra \Sel(\cG_{K\pr_{v\pr}}) \quad \quad \text{and} \quad\quad\Sel(\cH_{K_{v}}) \ra \Sel(\cH_{K\pr_{v\pr}}) \quad \quad \quad \quad \ \ 
\ee
whenever $v\pr \not\in U\pr$ and $v$ is the place below $v\pr$. We write 
\[
\Sel(\cG_{U\pr}) \qq \text{and} \qq \Sel(\cH_{U\pr}) \qqqq \ \, \,
\]
(resp.,~$\Sel(\cG)$ and $\Sel(\cH)$ if $K\pr = K$) for the Selmer groups that result by completing \eqref{Sel-Sel} to Selmer conditions as in \eqref{comp-cond}. These Selmer groups are finite due to \Cref{Sel-fin}. Due to \eqref{Sel-comp}, for each $K\pr \in \sS$ restriction maps induce the maps
\[
\ \, \Sel(\cG) \ra \Sel(\cG_{U\pr}) \quad \quad \text{and} \quad \quad \Sel(\cH) \ra \Sel(\cH_{U\pr}). \quad \quad \quad \quad\ \ \,
\]
As in \S\ref{OC}, shrinking $U$ affects neither the above setup, nor the Selmer groups.
\epp

\bthm \lab{Sel-gr}
In the setup of \S\ref{Sel-BC}, if $\cG_K$ is \'{e}tale, $p$ is a prime dividing $\#\cG_K$, and $\sS$ consists of the $\bZ/p\bZ$-subextensions of $\ov{K}/K$, then 
\[
\#\Sel(\cG_{U\pr})
\]
is unbounded when $K\pr$ ranges in $\sS$.
\ethm

\bpf
Let $V\subset U$ be a nonempty open. Initial segments of the exact sequences of \Cref{PT-gen} for $\cG_{V}$ and $\cG_{V\pr}$ fit into the commutative diagram
\[
\xymatrix@C=12pt @R=16pt{
0 \ar[r] & \Sel(\cG) \ar[r]\ar[d]^-{a} & H^1(V, \cG) \ar[r]^-{l}\ar[d]^-{b} & \p{\bigoplus_{v\not\in U} \f{H^1(K_v, \cG)}{\Sel(\cG_{K_v})}} \oplus \p{\bigoplus_{v\in U \setminus V} \f{H^1(K_v, \cG)}{H^1(\cO_v, \cG)}} \ar[d]^-{c} \ar[rr]^-{y(\cG_V)} & & \Sel(\cH)^* \\
0 \ar[r] & \Sel(\cG_{U\pr}) \ar[r] & H^1(V\pr, \cG) \ar[r]^-{l\pr} & \p{\bigoplus_{v\pr\not\in U\pr} \f{H^1(K\pr_{v\pr}, \cG)}{\Sel(\cG_{K\pr_{v\pr}})}} \oplus \p{\bigoplus_{v\pr\in U\pr \setminus V\pr} \f{H^1(K\pr_{v\pr}, \cG)}{H^1(\cO_{v\pr}, \cG)}} \ar[rr]^-{y(\cG_{V\pr})} & & \Sel(\cH_{U\pr})^*.
}
\]
\bcl \lab{b-bdd}
As $V$ and $K\pr$ vary, $\#\Ker b$ is bounded.
\ecl

\bpf
By the injectivity aspect of \Cref{Cart}, 
\[
H^1(V, \cG) \subset H^1(K, \cG) \qq \text{and}\qq H^1(V\pr, \cG) \subset H^1(K\pr, \cG),
\]
so $\Ker b \subset H^1(\Gal(K\pr/K), \cG(K\pr))$. It remains to observe that the cardinality of $H^1(\Gal(K\pr/K), \cG(K\pr))$ is bounded in terms of $p$ and $\#\cG$.
\epf

\bcl \lab{c-unbdd}
As $V$ and $K\pr$ vary, $\#\Ker c$ is unbounded.
\ecl

\bpf
We fix an $m \in \bZ_{> 0}$. Since $\cG$ is finite \'{e}tale over a nonempty open of $U$, \v{C}ebotarev density theorem gives a set $\Sigma$ of $m$ closed points $v \in U$ for which $\mu_p(K_v) = \mu_p(\ov{K}_v)$ and $\cG_{\cO_v}$ is constant. 

We fix a $v \in \Sigma$ and let $\underline{\bZ/p^r \bZ}_{\cO_v}$ be a direct summand of $\cG_{\cO_v}$. Since $H^1(K_v, \bZ/p^r\bZ)$ is the group of homomorphisms 
\[
h\colon \Gal(\ov{K}_v/K_v) \ra \bZ/p^r\bZ
\]
and $H^1(\cO_v, \bZ/p^r\bZ)$ is the subgroup of unramified $h$, every ramified $\bZ/p\bZ$-extension $K\pr_{v\pr}/K_v$ kills a nonzero element of $H^1(K_v, \cG)/H^1(\cO_v, \cG)$. We fix such a $K\pr_{v\pr}/K_v$: there are many to choose from if $\Char K_v = p$, and there is at least one if $\Char K_v \neq p$ due to the $\mu_p(K_v) = \mu_p(\ov{K}_v)$ requirement. 

We use \cite{NSW08}*{9.2.8} to find a $\bZ/p\bZ$-subextension $\ov{K}/K\pr/K$ that interpolates the chosen local extensions $K\pr_{v\pr}/K_v$ and set $V \ce U - \Sigma$ to arrive at a $c$ with 
\[
\#\Ker c \ge p^m. \qedhere
\]
\epf

Since $\#\Coker l$ is bounded by $\#\Sel(\cH)$, unboundedness of $\#\Ker c$ supplied by \Cref{c-unbdd} implies that of $\#\Ker(c|_{\im l})$. By \Cref{b-bdd}, $\#\Ker b$ stays bounded, so unboundedness of $\#\Ker(c|_{\im l})$ implies that of $\#\Coker a$, i.e., that of $\#\Sel(\cG_{U\pr})$ when $K\pr$ ranges in $\sS$.
\epf

\bcor \lab{CG-unbdd}
For a prime $p$, the cardinalities 
\[
\#\Pic(S\pr)[p]
\]
are unbounded when $K\pr$ ranges over the $\bZ/p\bZ$-extensions of $K$.
\ecor

\bpf
We use \cite{Ces15a}*{B.1 (a)} to replace $\#\Pic(S\pr)[p]$ by $\#H^1(S\pr, \bZ/p\bZ)$. Then it remains to apply \Cref{Sel-gr} to $U = S$ and $\cG = \bZ/p\bZ$ with $\Sel(\cG_{K\pr_{v\pr}}) = H^1(K\pr_{v\pr}, \bZ/p\bZ)$ for $v\pr \mid \infty$.
\epf

\brem
For further results similar to \Cref{CG-unbdd}, see, for instance, \cite{Mad72}.
\erem

\bcor \lab{phi-unbdd}
For an isogeny $\phi\colon A \ra B$ between abelian varieties over $K$ and a prime $p$ that divides the order of some $K$-\'{e}tale subgroup $G \subset A[\phi]$ (if $p \neq \Char K$, then the $K$-\'{e}taleness of $G$ is automatic),  
\[
\#\Sel_{\phi} A_{K\pr}
\]
is unbounded when $K\pr$ ranges over the $\bZ/p\bZ$-extensions of $K$.
\ecor

\bpf
Let $\psi\colon A \ra C$ be an isogeny with kernel $G$. Since 
\[
\#\Ker(\Sel_{\psi} A_{K\pr} \ra \Sel_{\phi} A_{K\pr})
\]
is bounded by $\#(A[\phi]/G)$, we may assume that $\psi = \phi$. In this case, we let $\phi^\vee$ be the dual isogeny and choose $U$ and \eqref{Sel-Sel} as in \Cref{CPT-phi} (using \Cref{boring-check}) to argue that \Cref{Sel-gr} applies.
\epf

In characteristic $p$, \Cref{phi-unbdd} may be supplemented by the following result.

\bthm \lab{p-Sel-unbdd}
For a prime $p$ and a nonzero abelian variety $A$ over $K$, there is an $n \in \bZ_{> 0}$ for which
\[
\#\Sel_{p^n} A_{K\pr}
\]
is unbounded when $K\pr$ ranges over the $\bZ/p^n\bZ$-extensions of $K$. Moreover, one may choose $n = 1$ if $A[p](\ov{K}) \neq 0$ (for instance, if $\Char K \neq p$) or if $A$ is supersingular (for instance, if $A[p](\ov{K}) = 0$ and $\dim A \le 2$).
\ethm

\brem
For every $N \ge n$, the sizes of the kernel and the cokernel of the map
\[
\qq \Sel_{p^{n}} A_{K'} \ra (\Sel_{p^{N}} A_{K'})[p^{n}]
\]
are bounded by $p^{2ng}$, see \cite{Ces15a}*{6.7 (a)}. Therefore, \Cref{p-Sel-unbdd} also shows that $\#\Sel_{p^N} A_{K\pr}$ is unbounded when $K'$ ranges over the $\bZ/p^{n}\bZ$-extensions of $K$.
\erem

In the case when $A[p](\ov{K}) = 0$, the proof of \Cref{p-Sel-unbdd} will use the following lemmas.

\blem \lab{slopes}
For a connected-connected $p$-divisible group $G$ over a field $k$ of characteristic $p$, there exist integers $x > y > 0$ and $z \ge 0$ such that
\[
\Ker(\Frob_{p^{xt},\, G/k}) \subset G[p^{yt + z}] \qq \text{for every} \quad t \in \bZ_{\ge 0}.
\]
\elem

\bpf
For any $x$, $y$, and $z$, the indicated inclusions may be tested over the algebraic closure of $k$, so we loose no generality by assuming that $k = \ov{k}$. Since $G$ is connected-connected, all its Dieudonn\'{e}--Manin slopes lie in the open interval $(0, 1)$.

If $G$ is isoclinic of slope $\f{r}{s}$, then its $s$-fold relative Frobenius morphism $\Frob_{p^{s},\, G/k}$ identifies with multiplication by $p^r$, so it suffices to set $x \ce s$, $y \ce r$, and $z \ce 0$ (with these choices the indicated inclusions are even equalities). Thus, more generally, if $G$ is a product of isoclinic $p$-divisible groups $G_i$ with slopes $\{\f{r_i}{s_i}\}$ and we set $x \ce \prod s_i$, then 
\[
\Ker(\Frob_{p^{xt},\, G_i/k}) = G_i[p^{r_i \cdot \f{x}{s_i} \cdot t}] \qq \text{for every $i$ and every $t \in \bZ_{\ge 0}$},
\]
so it suffices to in addition set $y \ce \max_i \{ r_i \cdot \prod_{i' \neq i} s_{i'} \}$ and $z \ce 0$.

In general, thanks to the Dieudonn\'{e}--Manin classification and the assumption $k = \ov{k}$, there is a $k$-isogeny $f\colon G \ra G'$ towards a $p$-divisible group $G'$ that is a product of isoclinic $p$-divisible groups. Thus, if $x$ and $y$ are chosen for $G'$ as in the previous paragraph, then, in order to obtain the sought triple $x$, $y$, $z$ for $G$, it remains to let $z$ be any nonnegative integer such that $\Ker f \subset G[p^z]$.
\epf

\blem \lab{ramify}
For every $n, \ell \in \bZ_{\ge 1}$, a nonarchimedean local field $k$ of characteristic $p > 0$ has infinitely many totally ramified $\bZ/p^n\bZ$-extensions $\wt{k}/k$ such that $\Gal(\wt{k}/k)$ acts trivially on $\wt{\fo}/\wt{\fm}^{\ell}$, where $\wt{\fo}$ denotes the ring of integers of $\wt{k}$ and $\wt{\fm} \subset \wt{\fo}$ denotes the maximal ideal.
\elem

\bpf
The condition on the triviality of the action means that the ramification subgroup $\Gal(\wt{k}/k)_{\ell - 1}$ equals the entire $\Gal(\wt{k}/k)$. In terms of the upper numbering, this means that 
\[
\tst \Gal(\wt{k}/k)^{\varphi_{\wt{k}/k}(\ell - 1)} = \Gal(\wt{k}/k), \qq \text{where } \quad \varphi_{\wt{k}/k}(\ell - 1) = \int_0^{\ell - 1} \f{dx}{[\Gal(\wt{k}/k)_0\, : \, \Gal(\wt{k}/k)_x]}.
\]
Thus, since $\varphi_{\wt{k}/k}(\ell - 1) \le \ell - 1$, it suffices to ensure that 
\[
\Gal(\wt{k}/k)^{\ell - 1} = \Gal(\wt{k}/k), \qq \text{or that even} \qq \Gal(\wt{k}/k)^{\ell} = \Gal(\wt{k}/k).
\]
Local class field theory then reduces us to finding infinitely many continuous $\bZ/p^n\bZ$-quotients of $\fo^\times$ onto which $1 + \fm^{\ell}$ maps surjectively, where $\fo$ denotes the ring of integers of $k$ and $\fm \subset \fo$ denotes the maximal ideal. By \cite{Neu99}*{II.5.7~(ii) and its proof}, there is a topological group isomorphism 
\[
\tst \fo^\times \simeq (\fo/\fm)^\times \times \prod_{i = 1}^\infty \bZ_p
\]
such that the image of the subgroup $1 + \fm^\ell \subset \fo^\times$ contains the subgroup $\prod_{i = 1}^{f(\ell)} 0 \times \prod_{i = f(\ell) + 1}^\infty \bZ_p$ for some $f(\ell) \in \bZ_{\ge 1}$. It remains to observe that the quotient $\bZ_p/p^n\bZ_p$ of the $i\th$ copy of $\bZ_p$ with $i > f(\ell)$ gives rise to a sought quotient of $\fo^\times$.
\epf

\bpf[Proof of \Cref{p-Sel-unbdd}]
\Cref{phi-unbdd} settles the case $\Char K \neq p$, so we assume that $\Char K = p$. 

The overall structure of the argument will be similar to the one used to prove \Cref{Sel-gr}. Namely, we let $V \subset U \subset S$ be nonempty opens such that $A$ extends to an abelian scheme $\cA \ra U$, and we use the sequence \eqref{Coates} with $\phi = [p^n]_A$ (and an $n$ to be fixed later) to obtain the commutative diagram 
\[
\xymatrix@C=12pt @R=16pt{
0 \ar[r] & \Sel_{p^n} A \ar[r]\ar[d]^-{a} & H^1(V, \cA[p^n]) \ar[r]^-{l}\ar[d]^-{b} &  \displaystyle \bigoplus_{v\not\in  V} H^1(K_v, A)[p^n] \ar[d]^-{c} \ar[r] & (\Sel_{p^n} A^\vee)^* \\
0 \ar[r] & \Sel_{p^n} A_{K\pr} \ar[r] & H^1(V\pr, \cA[p^n]) \ar[r]^-{l\pr} &  \displaystyle\bigoplus_{v\pr\not\in V\pr} H^1(K\pr_{v\pr}, A)[p^n]
}
\]
with exact rows.  As in the final paragraph of the proof of \Cref{Sel-gr}, due to the uniform boundedness of $\#\Ker b$ supplied by the proof of \Cref{b-bdd}, it suffices to prove that there is an $n$ subject to the constraints of the last sentence of the claim such that $\#\Ker c$ is unbounded when $V$ and $K\pr$ vary. Moreover, since we may use \cite{NSW08}*{9.2.8} to interpolate any finite set of local $\bZ/p^n\bZ$-extensions by a global $\bZ/p^n\bZ$-extension, it suffices to show that for some $n$ satisfying the constraints and for infinitely many $v \in U$ there is a $\bZ/p^n\bZ$-extension $K'_{v'}/K_v$ such that 
\be \lab{wanted}
\Ker(H^1(K_v, A)[p^n] \ra H^1(K\pr_{v\pr}, A)[p^n]) \neq 0.
\ee
By Tate local duality, the requirement \eqref{wanted} is equivalent to the requirement that the norm map
\be \lab{wanted-2}
\cN \colon A^\vee(K\pr_{v\pr}) \ra A^\vee(K_v) \qqq  \text{is not surjective.}
\ee
To find infinitely many $v\in U$ satisfying \eqref{wanted-2} for some $K'_{v'}$, we split the argument into cases.

\emph{The case when $A[p](\ov{K}) \neq 0$.} 
In this case, neither of the isogenous $p$-divisible groups $A[p^\infty]$ and $A^\vee[p^\infty]$ is connected, so $A^\vee[p](\ov{K}) \neq 0$, too. Therefore, at the cost of shrinking $U$, we may assume that $\cA^\vee[p]$ is an extension of a nonzero finite \'{e}tale $U$-group scheme $\cQ$ by a finite flat $U$-group scheme $\cH$ that has connected fibers. By the \v{C}ebotarev density theorem, there are infinitely many $v \in U$ such that $\cQ(\bF_v) \neq 0$. For such $v$, due to the vanishing of $H^1(\bF_v, \cH)$ supplied, for instance, by \cite{Ces15d}*{5.7~(b)}, one also has $\cA^\vee(\bF_v)[p] \neq 0$.

We choose $n = 1$, let $v$ be such that $\cA^\vee(\bF_v)[p] \neq 0$, and let $K'_{v'}/K_v$ be any ramified $\bZ/p\bZ$-extension. The norm map $\cN$ of \eqref{wanted-2} reduces to multiplication by $p$ on $\cA^\vee(\bF_v)$, so cannot be surjective because $\#(\cA^\vee(\bF_v)/p\cA^\vee(\bF_v)) = \#\cA^\vee(\bF_v)[p]$.

\emph{The case when $A$ is supersingular.}
We choose $n = 1$ and for any $v \in U$ use \Cref{ramify} to choose a ramified $\bZ/p\bZ$-extension $K'_{v'}/K_v$ for which $\Gal(K'_{v'}/K_v)$ acts trivially on $\cO_{v'}/\fm_{v'}^{p^2}$. For proving that the norm map $\cN$ of \eqref{wanted-2} is not surjective, we consider the formal group of $A^\vee_{K_v}$, i.e.,~the formal completion of $\cA^\vee_{\cO_v}$ along the identity section of the special fiber. The $\cO_v$-points of this formal group identify with
\[
\tst \varprojlim_{n \ge 1} \Ker(\cA^\vee(\cO_v/\fm_v^n) \surjects \cA^\vee(\bF_v)),
\]
and hence also with the kernel
\[
\Ker\p{\cA^\vee(\cO_v) \surjects \cA^\vee(\bF_v)}
\]
of the reduction map, and likewise over $\cO_{v'}$. We seek to show that
\be\lab{wanted-3}
\tst \cN\p{\Ker(\cA^\vee(\cO_{v'}) \surjects \cA^\vee(\bF_{v'}) )} \subset \Ker(\cA^\vee(\cO_v) \surjects \cA^\vee(\cO_v/\fm_v^p))
\ee
(the indicated surjectivity results from the $\cO_v$-smoothness of $\cA^\vee_{\cO_v}$). Once this is done, it will follow that $\cN\colon \cA^\vee(\cO_{v'}) \ra \cA^\vee(\cO_v)$ cannot be surjective because $\#\cA^\vee(\bF_{v'}) < \# \cA^\vee(\cO_v/\fm_v^p)$. 

Since $\cA^\vee(\cO_v/\fm_v^p) \subset \cA^\vee(\cO_{v'}/\fm_{v'}^{p^2})$, for \eqref{wanted-3} it suffices to show that $\cN$ induces the zero map on
\[
\Ker(\cA^\vee(\cO_{v'}/\fm_{v'}^{p^2}) \ra \cA^\vee(\bF_{v'})).
\]
This induced map agrees with multiplication by $p$ because $\Gal(K'_{v'}/K_v)$ acts trivially on $\cO_{v'}/\fm_{v'}^{p^2}$. Moreover, since $A$ is supersingular, so is $A^\vee$, and hence the $p$-divisible group $A^\vee[p^\infty]$ is isoclinic of slope $\f{1}{2}$. In particular, the multiplication by $p$ map of $A^\vee$ identifies with the relative $p^2$-Frobenius morphism $\Frob_{p^2,\, A^\vee/K}$. Therefore, the map induced by $\cN$ on $\Ker(\cA^\vee(\cO_{v'}/\fm_{v'}^{p^2}) \ra \cA^\vee(\bF_v))$ agrees with the map induced by the relative $p^2$-Frobenius of $\cA^\vee_{\cO_{v'}}$, and hence vanishes.

\emph{The case when $A[p](\ov{K}) = 0$.} 
In this case, the isogenous $p$-divisible groups $A[p^\infty]$ and $A^\vee[p^\infty]$ are connected-connected, so \Cref{slopes} provides integers $x > y > 0$ and $z \ge 0$ such that 
\be \lab{slope-input}
\Ker(\Frob_{p^{xt},\, A^\vee/k}) \subset A^\vee[p^{yt + z}] \qq \text{for every} \quad t \in \bZ_{\ge 0}.
\ee
We choose a $t$ for which $xt > yt + z$, set $n \ce yt + z$, and for a $v \in U$ use \Cref{ramify} to choose a totally ramified $\bZ/p^n\bZ$-extension $K'_{v'}/K_v$ for which $\Gal(K'_{v'}/K_v)$ acts trivially on $\cO_{v'}/\fm_{v'}^{p^{xt}}$. Similarly to the proof of the supersingular case, we seek to show that
\be \lab{wanted-4}
\tst \cN\p{\Ker(\cA^\vee(\cO_{v'}) \surjects \cA^\vee(\bF_{v'}) )} \subset \Ker(\cA^\vee(\cO_v) \surjects \cA^\vee(\cO_v/\fm_v^{p^{xt - n}})),
\ee
which will prove the sought nonsurjectivity of $\cN$ because $\#\cA^\vee(\bF_{v'}) < \# \cA^\vee(\cO_v/\fm_v^{p^{xt - n}})$. To prove \eqref{wanted-4}, it suffices to prove that the map induced by $\cN$ on 
\[
\Ker(\cA^\vee(\cO_{v'}/\fm_{v'}^{p^{xt}}) \ra \cA^\vee(\bF_{v'})).
\]
is zero. Since $\Gal(K'_{v'}/K_v)$ acts trivially on $\cO_{v'}/\fm_{v'}^{p^{xt}}$, this induced map is multiplication by $p^n$, so, due to \eqref{slope-input}, it factors through the zero map induced by the relative $p^{xt}$-Frobenius of $\cA^\vee_{\cO_{v'}}$.
\epf

\brem
Since the sequence 
\[
0 \ra (\Sel_{p^n} A_{K\pr})^{\Gal(K\pr/K)} \ra  H^1(V\pr, \cA[p^n])^{\Gal(K\pr/K)} \xra{l\pr} \p{\im l\pr}^{\Gal(K\pr/K)}
\]
is exact, the proof of \Cref{p-Sel-unbdd} shows that even $\#(\Sel_{p^n} A_{K\pr})^{\Gal(K\pr/K)}$ is unbounded (and similarly in \Cref{phi-unbdd}). 
\erem


\begin{appendix}

\section{Continuity of cup products} \lab{app-a}

The goal of this appendix is to prove that cup product pairings on local cohomology groups are continuous, see \Cref{cup-cont} for a precise statement. Such continuity is implied by the assertion \cite{Mil06}*{III.6.5 (e)} (whose proof is omitted in loc.~cit.)~and is crucial for this paper through its roles in the Tate--Shatz local duality \cite{Mil06}*{III.6.10} and in the proof of \Cref{im-closed}.

As always, the topology on cohomology is that defined in \cite{Ces15d}*{3.1--3.2}. However, we also use \cite{Ces15d}*{5.11 and 6.5 (with 3.5 (d))}, which guarantee agreement with the ``\v{C}ech topology.'' In \Cref{H1-c-i}, we recall the needed \v{C}ech-theoretic notation; see \cite{Ces15d}*{5.1} for further recollections.

\blem \lab{H1-nbhd}
For a local field $k$, a $k$-group scheme $G$ locally of finite type, and an $x \in H^1(k, G)$, there is a finite extension $k\pr/k$ such that the image of $H^1(k\pr/k, G) \hra H^1(k, G)$ (consisting of classes of right $G$-torsors that become trivial over $k\pr$) contains an open neighborhood of $x$.
\elem

\bpf
By \cite{SGA3Inew}*{\upshape{VII}$_{\text{\upshape{A}}}$, 8.3}, $G$ is an extension 
\[
1 \ra H \ra G \ra Q \ra 1
\]
of a smooth $k$-group scheme $Q$ by a finite connected $H$. By \cite{Ces15d}*{3.5 (a) and 4.2}, the map
\[
H^1(k, G) \ra H^1(k, Q)
\]
has open fibers, so a $k\pr/k$ killing the fiber containing $x$ would suffice. To arrive at such a $k\pr$, we replace $k$ by a finite extension to kill the image of $x$ in $H^1(k, Q)$ and apply \cite{Ces15d}*{5.7~(b)}, which supplies a finite extension of $k$ that kills the entire $H^1(k, H)$.
\epf

\blem \lab{H1-c-i}
For a finite extension $k\pr/k$ of local fields and a finite $k$-group scheme $G$, the map 
\[
H^1(k\pr/k, G) \hra H^1(k, G)
\]
is a closed embedding (as in \cite{Ces15d}*{5.1}, we let $Z^1_{k\pr/k, G}$ be the $k$-scheme of $1$-cocycles and endow $H^1(k\pr/k, G)$ with the quotient topology via $Z^1_{k\pr/k, G}(k) \surjects H^1(k\pr/k, G)$).
\elem

\bpf
We fix an algebraic closure $\ov{k}$ containing $k\pr$, so 
\[
\tst H^1(k, G) = \varinjlim_{\ov{k}/\wt{k}/k\pr} H^1(\wt{k}/k, G),
\] where $\wt{k}$ ranges over the indicated finite subextensions. By \cite{Ces15d}*{5.11}, if each $H^1(\wt{k}/k, G)$ is topologized analogously to $H^1(k\pr/k, G)$, then the topology on $H^1(k, G)$ agrees with the direct limit topology. It therefore suffices to show, as we do below, that each 
\[
H^1(k\pr/k, G) \hra H^1(\wt{k}/k, G)
\]
is closed.

For $n \ge 0$, we set $\wt{k}_n \ce \tensor_{i = 0}^n \wt{k}$ (tensor product over $k$) and let 
\[
C^n_{\wt{k}/k, G} \ce \Res_{\wt{k}_n/k} (G_{\wt{k}_n})
\]
be the scheme of $n$-cochains (for $G$ with respect to $\wt{k}/k$). Since $G$ is finite, $C^n_{\wt{k}/k, G}$ is an affine $k$-group scheme of finite type. Thus, 
\[
C^1_{k\pr/k, G} \hra C^1_{\wt{k}/k, G}
\]
is a closed immersion by \cite{SGA3Inew}*{\upshape{VI}$_{\text{\upshape{B}}}$, 1.4.2}, so 
\[
Z^1_{k\pr/k, G} \hra Z^1_{\wt{k}/k, G}
\]
is one, too. It remains to note that the $C^0_{\wt{k}/k, G}(k)$-orbit quotient map 
\[
Z^1_{\wt{k}/k, G}(k) \surjects H^1(\wt{k}/k, G)
\]
is closed because $C^0_{\wt{k}/k, G}(k)$ is finite.
\epf

\bprop \lab{cup-cont}
For a local field $k$ and a bilinear pairing $G \times_k H \ra F$ of commutative $k$-group schemes locally of finite type with $G$ and $H$ finite, the cup product induces a continuous map 
\[
H^n(k, G) \times H^m(k, H) \ra H^{n + m}(k, F) \quad \quad \text{ for every } \quad \quad n, m \in \bZ_{\ge 0}.
\]
\eprop

\bpf
By \cite{Sha72}*{p.~208, Thm.~42}, every element of $H^n(k, G)$ lies in the image of $H^n(k\pr/k, G)$ for some finite extension $k\pr/k$, and likewise for $H^m(k, H)$. Moreover, $H^n(k, G)$ and $H^m(k, H)$ are discrete except for $H^1$. Working in neighborhoods of fixed elements of $H^n(k, G)$ and $H^m(k, H)$ and using \Cref{H1-nbhd,H1-c-i}, we therefore reduce to proving the continuity of the composition
\[
H^n(k\pr/k, G) \times H^m(k\pr/k, H) \xra{- \cup -} H^{n + m}(k\pr/k, F) \xra{y} H^{n + m}(k, F) \quad\quad \text{for every $k\pr/k$.}
\]
Continuity of $y$ is part of the agreement with the \v{C}ech topology, whereas the cup product lifts to a map
\[
Z^n_{k\pr/k, G}(k) \times Z^m_{k\pr/k, H}(k) \ra Z^{n + m}_{k\pr/k, F}(k)
\]
that is continuous because it is induced by a $k$-scheme morphism.
\epf


\section{Selmer conditions for dual isogenies are orthogonal complements} \lab{app-b}

We seek to justify the legitimacy of the choice of Selmer conditions in Example \ref{CPT-phi} by proving \Cref{boring-check}, which is standard but seems to lack a reference.

\bprop \lab{boring-check}
For a local field $k$ and dual isogenies $\phi \colon A \ra B$ and $\phi^\vee \colon B^\vee \ra A^\vee$ between abelian varieties over $k$, the subgroups 
\[
B(k)/\phi A(k) \subset H^1(k, A[\phi]) \quad \quad \text{and} \quad \quad A^\vee(k)/\phi^\vee B^\vee(k) \subset H^1(k, B^\vee[\phi^\vee])
\] 
are orthogonal complements under the Tate-Shatz local duality pairing 
\[
H^1(k, A[\phi]) \times H^1(k, B^\vee[\phi^\vee]) \ra \bQ/\bZ.
\]
\eprop

\bpf
Due to the commutativity of diagrams such as the first one in \cite{Mil06}*{III.7.8}, the case $\phi = n_A$ for $n \in \bZ_{> 0}$ is an implicit corollary of the proofs \cite{Mil06}*{I.3.4, I.3.7, and III.7.8} of Tate local duality for abelian varieties. Thus, we assume the $\phi = n_A$ case to be known and deduce the general case.

By symmetry, it suffices to show that $B(k)/\phi A(k)$ is the annihilator of $A^\vee(k)/\phi^\vee B^\vee(k)$. We set $n \ce \deg \phi$ and let $\psi \colon B \ra A$ be the isogeny for which $\psi \circ \phi = n_A$. In the diagram
\vspace{0.5em} 
\be\ba\lab{curveball}
\xymatrix{
B(k)/\phi A(k) \ar[d]^{\psi} \ar@{^(->}[r]  & H^1(k, A[\phi]) \ar[d]  \ar@{}[r]|-{\bigtimes} & H^1(k, B^\vee[\phi^\vee]) \ar@/^1pc/[rr] & A^\vee(k)/\phi^\vee B^\vee(k) \ar@{_(->}[l] & \bQ/\bZ \ar@{=}[d] \\
A(k)/nA(k) \ar@{^(->}[r]  & H^1(k, A[n]) \ar@{}[r]|-{\bigtimes} & H^1(k, A^\vee[n]) \ar[u]_-{H^1(k, \psi^\vee)} \ar@/_1pc/[rr] & A^\vee(k)/nA^\vee(k) \ar@{_(->}[l] \ar@{->>}[u] & \bQ/\bZ
}
\ea\ee

in which $\psi^\vee$ is the isogeny dual to $\psi$, the curved arrows are the Tate--Shatz local duality cup product pairings, and the wing squares commute. By \cite{Oda69}*{1.1}, the inclusion $A[\phi] \hra A[n]$ identifies with the Cartier dual of $\psi^\vee\colon A^\vee[n] \ra B^\vee[\phi^\vee]$, so the commutativity of the pairing square results from using \cite{GH71}*{3.1} to identify cup product pairings with Ext-product pairings, using \cite{GH70}*{4.5} to identify Ext-product pairings with Yoneda edge product pairings, and using the commutativity of 
\[
\xymatrix{
 \Ext^1(B^\vee[\phi^\vee], \bG_m) \ar[d]_{\Ext^1(\psi^\vee, \bG_m)}  \ar@{}[r]|-{\bigtimes} & \Ext^1(\bZ, B^\vee[\phi^\vee]) \ar[r] &  \Ext^2(\bZ, \bG_m) \ar@{=}[d] \\
 \Ext^1(A^\vee[n], \bG_m) \ar@{}[r]|-{\bigtimes} & \Ext^1(\bZ, A^\vee[n]) \ar[r] \ar[u]^{\Ext^1(\bZ, \psi^\vee)}  &  \Ext^2(\bZ, \bG_m)
}
\]
that results from interpreting $\Ext$'s as $\Hom$'s in a derived category.

The commutativity of \eqref{curveball} and the assumed $\phi = n_A$ case show that $B(k)/\phi A(k)$ kills $A^\vee(k)/\phi^\vee B^\vee(k)$. Moreover, an $x\in H^1(k, A[\phi])$ that kills $A^\vee(k)/\phi^\vee B^\vee(k)$ maps to $A(k)/nA(k)$ in $H^1(k, A[n])$, so $x \in B(k)/\phi A(k)$ due to the commutativity of the following diagram with exact rows:
\[
\begin{gathered}[b]
\xymatrix{
B(k)/\phi A(k) \ar[d]^{\psi} \ar@{^(->}[r]  & H^1(k, A[\phi]) \ar[d]  \ar@{->>}[r] & H^1(k, A)[\phi] \ar@{^(->}[d] \\
A(k)/nA(k) \ar@{^(->}[r]  & H^1(k, A[n]) \ar@{->>}[r] & H^1(k, A)[n]. 
}
\\[-\dp\strutbox]
\end{gathered}
\qedhere
\]
\epf

\end{appendix}

\end{document}